\documentclass{tac}

\usepackage[all]{xy}

\usepackage{amssymb}

\usepackage[colorlinks=true]{hyperref}

\author{Panagis Karazeris and Konstantinos Tsamis}

\address{
 Department of Mathematics,
University of Patras\\ Patras, Greece\\
}

\title {Regular and effective regular categories of locales}

\copyrightyear{2020}


\eaddress{pkarazer@upatras.gr \CR tsamiskonstantinos@gmail.com }

\newtheorem{theorem}{Theorem}

\newtheoremrm{rem}{Remark}

\let\pf\proof
\let\epf\endproof
\mathrmdef{Hom} \mathbfdef{Set}

\begin{document}
\maketitle

\begin{abstract}
We examine the analogues for the respective categories of locales of two well-known results about regularity and effectiveness of some categories of spaces. We show that the category of compact regular locales is effective regular (=Barr-exact). We also show that the category of compactly generated Hausdorff locales is regular, provided that it is coreflective within Hausdorff locales. We do not appeal to the existence of points (which would render the two results trivial) but rely on the treatment of the subject by methods that are valid in the internal logic of a topos. On the course to the result about compactly generated locales we arrive at a generalization of a result of B. Day and R. Street, deriving regularity for a cocomplete category containing a dense regular subcategory closed under finite limits and colimits and satisfying a certain compatibility condition of pullbacks  with appropriate colimits.
\end{abstract}


\section{Introduction}
While regular and, even more so, effective regular categories occur more frequently in the realm of algebra there are two well-known cases of categories of spaces that have these features. The category of compact Hausdorff spaces is effective regular and the category of compactly generated (weakly) Hausdorff spaces is regular \cite{cmv}. It is a rather natural question to ask whether the corresponding categories of locales maintain these features. 

For the case of compact Hausdorff locales we know from \cite{t} that they form a regular category. We show here that it is also effective. The extra step, effectiveness of equivalence relations, almost exists implicitly in the work of \cite{v} on proper maps of locales, in particular his result that proper closed equivalence relations on compact locales are effective. 

The situation concerning compactly generated Hausdorff locales is much more complicated. First of all we adopt the definition of compactly generated locales introduced in \cite{e}, which constitutes the main, if not the only, study of such locales: a Hausdorff locale is compactly generated if it is isomorphic to the colimit of the (directed, extremal monomorphic) diagram of its compact sublocales (via the canonical comparison as a co-cone for that diagram). The major question that is left open in that work is whether such locales form a coreflective subcategory of that of Hausdorff locales. This  would be the case if, for every Hausdorff locale, the colimit in question was Hausdorff (as for example would be the case if the canonical comparison described above was monomorphic), in which case the comparison map would be the co-unit of the adjunction. The question of coreflectivity is important for the way products (and hence also pullbacks) are calculated in that category, namely whether they are calculated by applying coreflection to the localic product. This in turn affects our strategy for approaching the question of regularity of compactly generated Hausdorff locales. For that we adapt the argument due to \cite{ds} for deriving regularity of the inductive completion of a category from the regularity of the given category. The argument is familiar in the theory of locally presentable categories but its essential ingredients do not require local presentability. We present a generalization of the relevant result in \cite{ds} showing regularity for a a cocomplete category containing a dense regular subcategory closed under finite limits and colimits and satisfying a certain compatibility condition of pullbacks  with appropriate colimits. One key step is the ``uniformity lemma", namely that if objects (like the compactly generated locales) are built up from building blocks (like their compact sublocales), then every finite diagram of such objects can be expressed as a colimit of diagrams of the same type with the vertices among the building blocks (and all objects of the given diagram to be presented as colimits of building blocks over the same indexing category). This uses only the density of the building blocks and their closure under finite colimits in the broader category. Another property used in their proof, in connection with the existence of regular epi - mono factorizations and the stability of regular epis under pullback, is the commutation of pullbacks with a particular type of colimits (directed extremal monomorphic ones, in our case). We show that it is sufficient that the canonical comparison between the colimit of pullbacks to the pullback of the colimits is epimorphic. This is where the nature of the products plays a role. If we assume coreflectivity we arrive at that epimorphicity result and subsequently at the regularity of compactly generated Hausdorff locales. 

Our terminology is, we believe, standard. A map of locales $f \colon X \to Y$ is determined by a map $f^* \colon OY \to OX$ between the respective frames that preserves finite infima and all suprema. Hence it has a right adjoint $f_* \vdash f^* .$ The map is a surjection if $f^*$ reflects order. It is proper  if $f_*$ preserves directed suprema and, for all $U \in OX,$ $V \in OY ,$ $f_* (U \vee f^* V)= f_* U \vee V .$ Under the equivalence of the category of locales over $X$ with that of locales internal in sheaves on $X,$ proper maps in the former correspond to compact locales in the latter \cite{el}. A locale  $X$ is Hausdorff if its diagonal $X \to X\times X$ is closed. The category of Hausdorff locales is a reflective subcategory of the category of locales. The reflection associates with a given locale the one that is determined by the largest Hausdorff subframe of the underlying frame of the given locale (\cite{ch}, Theorem 1.2.2).

\paragraph{Acknowledgements} The second named author benefited from the co-financed by Greece and the European Union (European Social Fund- ESF) Operational Programme ``Human Resources Development, Education and Lifelong Learning'' in the context of the project ``Strengthening Human Resources Research Potential via Doctorate Research'' (MIS-5000432). Both authors are indebted to Vassilis Aravantinos - Sotiropoulos and Christina Vasilakopoulou for useful discussions. We are particularly thankful to the referee for pointing out an important mistake in an early version of this work, which resulted to a substantial revision of our initial approach, and for many other suggestions that improved the presentation of this work.

\section{Compactly generated Hausdorff locales}

Let us recall from \cite{e} that a Hausdorff locale $X$ (i.e one whose diagonal is closed) is called compactly generated if the canonical comparison map 
$$\varepsilon _X  \colon \mathrm{colim} C_i \to X, $$ 
where the colimit is taken over all the compact sublocales of $X$ (hence is a directed diagram of inclusions), is an isomorphism. 
For an arbitrary Hausdorff locale $X$ the above map is not known to be a monomorphism in the category of locales. In case it is, or at least the resulting colimit is Hausdorff, Escardo shows that it constitutes the counit of an adjunction, rendering the category $\mathrm{CGHLoc}$ of compactly generated Hausdorff locales a coreflective subcategory of Hausdorff locales. Let us refer to the assumption that $\mathrm{CGHLoc}$ is coreflective in Hausdorff locales as the  \textit{coreflectivity hypothesis}. We have 

\begin{proposition}
Under the coreflectivity hypothesis, if $(t_{ij} \colon X_i \to X_j )$ is the directed diagram of inclusions of the compact sublocales of a compactly generated Hausdorff one and $Y$ is any other such locale, then the canonical map 
$$ \mathrm{colim}_i ( X_i \times Y )  \to (\mathrm{colim}_i X_i  ) \times Y  $$ in $\mathrm{CGHLoc}$ is a split  epimorphism.  
\end{proposition}

\pf 
One has to be aware of the fact that even under the coreflectivity hypothesis a directed colimit of inclusions of Hausdorff locales that is calculated in the category of locales need not be Hausdorff, while a product of compactly generated locales, calculated in the category of locales, need not be compactly generated. The colimit $\mathrm{colim}_i X_i ,$ calculated in the category of locales, is by assumption compactly generated Hausdorff hence it is a colimit in the sense of $\mathrm{CGHLoc}.$ Its product with the Hausdorff locale $Y$ remains Hausdorff. The products $X_i \times Y,$ again calculated in the category of locales, have a factor which is compact Hausdorff. Compact Hausdorff locales are locally compact, hence exponentiable. Taking product with such a locale preserves colimits and we conclude that the products $X_i \times Y$ are compactly generated Hausdorff themselves (so they are products in $\mathrm{CGHLoc}$). On the other hand the colimit that occurs as the domain of the canonical morphism, does not coincide with the colimit taken in the category of locales since the latter need not be a Hausdorff locale. It is though an epimorphic image of the latter as explained at the end of the previous section. 

The product in the right hand side has to be the one in $\mathrm{CGHLoc}$ which means that we have (under the coreflectivity hypothesis) to apply the coreflection functor to the ordinary localic product. So it is isomorphic to the colimit $\mathrm{colim} _k C_k $ of the system of all compact sublocales of the localic product $(\mathrm{colim}_i X_i  ) \times Y .$ 
The inclusion of each $C_k$ into the localic product followed by the projection to $\mathrm{colim}_i X_i$ factors through a compact sublocale of that colimit, in particular $C_k  \rightarrow X_{i(k)} \rightarrowtail \mathrm{colim}_i X_i $ and similarly for the projection to the other factor, $C_k \rightarrow D_k \rightarrowtail Y.$ Hence the compact sublocales $X_{i(k)} \times D_k $ are final among the $C_k$ so their colimit is isomorphic to $(\mathrm{colim}_i X_i  ) \times Y $ in $\mathrm{CGHLoc}.$ On the other hand the maps $X_{i(k)} \times D_k \to X_{i(k)} \times Y \to \mathrm{colim}_i ( X_i \times Y )$ induce a map $ \mathrm{colim}_{k} (X_{i(k)} \times D_k ) \to \mathrm{colim}_i ( X_i \times Y )$ which fits in 
a commutative diagram (where we denote emphatically the colimit and product in $\mathrm{CGHLoc}$ where necessary).
$$\xy \xymatrix{ \ar[d]^{}
\mathrm{colim}_{k} (X_{i(k)} \times D_k )  \ar[r]^{\cong} & \mathrm{colim} _k C_k  \ar[r]^{\cong} & (\mathrm{colim}_i X_i  ) \times ^{CGH} Y \\
 \mathrm{colim}_i ( X_i \times Y ) \ar@{->>}[r]^{} & \mathrm{colim}_i^{GCH} ( X_i \times Y ) \ar[ur]^{}  & }
\endxy$$
where the horizontal composite is an isomorphism hence the canonical $ \mathrm{colim}_i ( X_i \times Y )  \to (\mathrm{colim}_i X_i  ) \times Y $ is a split epimorphism. 

\epf

\begin{proposition}
Under the coreflectivity hypothesis, if $Z$ is an object in $\mathrm{CGHLoc},$  $(t_{ij} \colon f_i \to f_j )$ is a directed diagram of inclusions between maps ($f_i \colon X_i \to Z$) with compact Hausdorff domain over it and $g\colon Y \to Z$ another map in the same category, then the canonical map 
$$ \mathrm{colim}_i ( X_i \times _Z Y )  \to (\mathrm{colim}_i X_i  ) \times _Z Y   $$ is a split epimorphism.  
\end{proposition}

\pf 
We want to apply the previous result relativized over a base $Z$, that is to exploit the previous result as a statement about products in $\mathrm{CGHLoc(Shv}Z).$ In order to do that, we need to make sure that the data of this Proposition, in particular the map $f \colon \mathrm{colim}_i X_i \to Z,$ give 
data of the previous Proposition when relativized over $Z,$ more specifically that $f \colon \mathrm{colim}_i X_i \to Z$ corresponds to a compactly generated locale in $\mathrm{Shv}(Z).$ Since this map is the colimit in $\mathrm{Loc}/Z$ of the $f_i \colon X_i \to Z ,$  we want to show that each such is a proper map. Indeed, for each for each $i$ the composite $X_i \to \mathrm{colim}_i X_i \to Z$ has a  factorization $X_i \to K_i  \rightarrowtail Z,$ where $K_i$ is a compact sublocale of $Z,$ which is proper: The map $X_{i} \to K_i$ is proper being a map between compact Hausdorff locales (\cite{t} 3.6.1), while $K_i  \rightarrowtail Z$ is proper being a closed inclusion (the image $K_i$ is closed as a compact sublocale of a Hausdorff one). Moreover, as a locale in $\mathrm{Shv}Z ,$ $X \to Z$ is Hausdorff when $X$ is because the diagonal in $\mathrm{Shv}Z ,$ $X \to X \times _Z X ,$ is closed when $X \to X \times X$ is. Obviously then the colimit of these composites $\mathrm{colim}_i X_i \to Z$ corresponds to a compactly generated locale in  $\mathrm{Loc(Shv}Z)$ and so does the map $Y\to Z$ by the same argument, hence we can apply the previous Proposition. Referring to the locales in $\mathrm{Shv}Z$ by the names of their corresponding maps to $Z,$ we notice that since the $f_i$ are compact, the localic products $f_i \times g $ in $\mathrm{Shv}Z$ are also products in $\mathrm{CGHLoc(Shv}Z).$ As explained earlier, their colimit in $\mathrm{CGHLoc(Shv}Z)$ will be an epimorphic image of the colimit in $\mathrm{Loc(Shv}Z).$ Seen as map in $\mathrm{Loc}/Z ,$ the domain of the latter will be just $\mathrm{colim}_i ( X_i \times _Z Y ) ,$ where the colimit is calculated in the category of locales. Both the domain of $\mathrm{colim}_i (f_i \times g)$ in $\mathrm{CGHLoc(Shv}Z),$ as well as $\mathrm{colim}_i ( X_i \times _Z Y )$ in $\mathrm{CGHLoc}$ arise as epimorphic images of that object.  The former maps epimorphically to the latter. This is so because the frame corresponding to it is a subframe of the one corresponding to the latter, both being subframes of $O(\mathrm{colim}_i ( X_i \times _Z Y )).$ This in turn is due to the fact that a subframe of a frame that happens to be under $OZ$ with the property that it has a closed diagonal in $\mathrm{Loc}$ will also be the (underlying frame of a) domain of a map whose diagonal is closed in $\mathrm{Loc(Sh}Z).$

Now by the previous proposition there exists a split epimorphism 
$$ \varepsilon \colon  \mathrm{colim}_i ^{CGH(ShZ)}(f_i \times g) \twoheadrightarrow \mathrm{colim}_i f_i \times ^{CGH(ShZ)} g $$
in $\mathrm{CGHLoc(Shv}Z),$ where the latter product is meant in the sense of this category. 
The domain of the map  $\mathrm{colim}_i f_i \times g $ as a product in $\mathrm{CGHLoc(Shv}Z)$ need not coincide with the product $\mathrm{colim}_i X_i \times _Z Y $ in $\mathrm{CGHaus}$ (there may exist sublocales $S\to Z$ of the locale $\mathrm{colim}_i X_i \times _Z Y  \to Z$ which are proper as maps to $Z$ but need not have a compact domain). But since $\mathrm{colim}_i X_i \times _Z Y$ is a cone for the discrete diagram formed by $\mathrm{colim}_i f_i$ and $g$ over $Z,$ there is a factorization through the (domain of the) product 
$$\beta \colon \mathrm{colim}_i X_i \times _Z Y \to \partial _0 (\mathrm{colim}_i f_i \times g) $$
 over $Z.$ It is easily seen that this is a monomorphism.
We have the following commutative diagram
$$\xy \xymatrix{  & \partial _0 (\mathrm{colim}_i ^{CGH(ShZ)}(f_i \times g))  \ar@{->>}[dd]^{\alpha} \ar@{->>}[r]^{\varepsilon} & \partial _0 (\mathrm{colim}_i f_i \times ^{CGH(ShZ)} g) \ar@/^1pc/[l]^{\mu}   \\
 \mathrm{colim}_i ( X_i \times _Z Y ) \ar@{->>}[ur]^{} \ar@{->>}[dr]^{} &  & \\
& \mathrm{colim}_i ^{CGH}( X_i \times _Z Y ) \ar[r]^{\gamma} & \mathrm{colim}_i ( X_i \times _Z ^{GCH}Y ) \ar@{>->}[uu]^{\beta}  & }
\endxy$$ 
Finally we get that, if $\mu$ is a splitting for $\varepsilon$ then $\alpha \cdot \mu \cdot \beta $ is a splitting for the comparison $\gamma \colon  \mathrm{colim}_i ( X_i \times _Z Y )  \to (\mathrm{colim}_i X_i  ) \times _Z Y .$
 \epf

\section{Regularity of the category of compactly generated Hausdorff locales}

We begin by generalizing a lemma due to B.~Day and R.~Street that is well-known for the case of locally presentable categories \cite{ds} . Its statement has only to do with density assumptions (of the presentable objects in the original case) and the closure of the dense subcategory under certain colimits. We include the proof for the sake of completeness of exposition. 

\begin{lemma}
Let $\mathcal{K}$ be a cocomplete category containing a dense subcategory $\mathcal{C}$ which is closed in $\mathcal{K}$ under finite colimits. Then for any small category with finite hom-sets $\mathcal{D}$ and diagram $D \in [\mathcal{D}, \mathcal{K}]$ we have that 
$$D \cong \mathrm{colim} \; ([\mathcal{D}, \mathcal{C}]\downarrow D \to \; [\mathcal{D}, \mathcal{C}] \to \; [\mathcal{D}, \mathcal{K}] )$$
\end{lemma}

\pf 
We show that, for all $d \in \mathcal{D},$ the evaluation at $d$ of the canonical morphism from the colimit to $D$ is an isomorphism in  
$\mathcal{K}.$ Colimits in $[\mathcal{D}, \mathcal{K}]$ are given object-wise so, denoting $i \colon \mathcal{C} \to \mathcal{K}$ the inclusion and $\partial _0 \colon [\mathcal{D}, \mathcal{C}]\downarrow D \to \; [\mathcal{D}, \mathcal{C}]$ the domain functor
\begin{eqnarray*}
\mathrm{colim} \; ([\mathcal{D}, \mathcal{C}]\downarrow D \to \; [\mathcal{D}, \mathcal{C}] \to \; [\mathcal{D}, \mathcal{K}] )(d)  &\cong & 
 \\
 ev_d (\; \mathrm{colim} \; ([\mathcal{D}, i]\cdot \partial _0 \; \colon \; [\mathcal{D}, \mathcal{C}]\downarrow D \to  [\mathcal{D}, \mathcal{C}] \to  [\mathcal{D}, \mathcal{K}] )\; )
 &\cong &  \\  
 \mathrm{colim}  ( \; ev_d \cdot [\mathcal{D}, i]\cdot \partial _0 \; \colon \; [\mathcal{D}, \mathcal{C}]\downarrow D \to \; [\mathcal{D}, \mathcal{C}] \to \; [\mathcal{D}, \mathcal{K}] \to \mathcal{K} ) & (1) &  \\
\end{eqnarray*}
On the other hand the density of $\mathcal{C}$ in $\mathcal{K}$ means that for all $d$ 
$$Dd \cong  \mathrm{colim} (\mathcal{C} \downarrow Dd \to \mathcal{C} \to \mathcal{K}),$$
while inspection gives that the composite 
$$ev_d \cdot [\mathcal{D}, i]\cdot \partial _0 \; \colon \; [\mathcal{D}, \mathcal{C}]\downarrow D \to \; [\mathcal{D}, \mathcal{C}] \to \; [\mathcal{D}, \mathcal{K}]  \to \mathcal{K} \quad (2) $$
is naturally isomorphic to the composite 
$$ i \cdot \partial _0 \cdot (ev_d \downarrow D) \; \colon \; [\mathcal{D}, \mathcal{C}]\downarrow D  \to \mathcal{C} \downarrow Dd \to \mathcal{C} \to \mathcal{K} \quad (3)$$
Moreover the functor $ev_d \downarrow D \colon  [\mathcal{D}, \mathcal{C}]\downarrow D \to \mathcal{C} \downarrow Dd $ is final, essentially because $ev_d \colon [\mathcal{D}, \mathcal{C}] \to \mathcal{C}$ has a left adjoint given by $C \mapsto \bigsqcup _{\mathcal{D}(d,-)} C$ (whose existence is granted by the fact that $\mathcal{D}$ has finite hom-sets and $\mathcal{C}$ is closed under finite colimits in $\mathcal{K}.$) Hence combining the isomorphisms $(1), \; (2), \; (3)$ with the latter finality result we get the desired isomorphism.
\epf 

The following proposition generalizes the main result of \cite{ds}. It relies on the possibility of conveniently writing a pullback diagram where the one leg is regular epi as colimit of diagrams with vertices in the dense subcategory. This is a key step for arriving at the characterization of regular locally finitely presentable categories in \cite{cpr} (cf. Lemma 12 and Theorem 14 there). We elaborate on it, showing that the previous lemma, using only density assumptions, suffices.

\begin{proposition}
Let $\mathcal{K} $ be a cocomplete and finitely complete category, such that it contains a dense regular subcategory $\mathcal{C}$ which is closed in $\mathcal{K}$ under finite limits and finite colimits. Assume that the objects of $\mathcal{K}$ are expressed (by density of $\mathcal{C}$) as colimits of objects from $\mathcal{C}$ of such kind that the canonical comparison map to the pullback of such colimits from the colimit of the pullbacks  of its components in $\mathcal{K}$ is a regular epimorphism. Then $\mathcal{K}$ is also regular.
\end{proposition}

\pf 
First we show the existence of regular epi- mono factorizations in $\mathcal{K}.$ We apply the above lemma for $\mathcal{D}$ the category $\bullet \to \bullet$ so that we express every morphism $X\to Y$ in $\mathcal{K}$ as a colimit of morphisms $X_i \to Y_i $ between objects in the full subcategory $\mathcal{C}.$ Using the regularity of $\mathcal{C}$ we take the regular epi - mono factorization $X_i \to W_i \to Y_i $ of every such morphism. Taking colimit of the appropriate kind we get a factorization 
$$X \cong \mathrm{colim}_i X_i \to  \mathrm{colim}_i W_i \to \mathrm{colim}_i Y_i\cong Y,$$
where the first morphism is obviously regular epi. We claim that the second one is mono because its kernel-pair consists of equal legs: Considering the pullbacks $W_i \times _{Y_i } W_i$ we have that the two legs to $W_i$ are equal since $W_ i \to Y_i$ is mono. Taking colimit over $i$ we get that the two outer morphisms $\mathrm{colim}_i (W_i \times _{Y_i } W_i ) \to \mathrm{colim}_i W_i$ in 

$$\xy \xymatrix{
\mathrm{colim}_i (W_i \times _{Y_i } W_i ) \ar[dr]^{} \ar[ddr]^{} \ar[drrr]^{}&   & & \\
   & \mathrm{colim}_i W_i  \times _{\mathrm{colim}_i Y_i } \mathrm{colim}_i W_i  \ar[rr]^{} \ar[d]^{} &  & \mathrm{colim}_i W_i  \ar[d]^{}\\
   & \mathrm{colim}_i W_i  \ar[rr]^{}  &   &  \mathrm{colim}_i Y_i }
\endxy $$
are equal. Hence the two composites 
$$\mathrm{colim}_i (W_i \times _{Y_i } W_i ) \to \mathrm{colim}_i W_i  \times _{\mathrm{colim}_i Y_i } \mathrm{colim}_i W_i  \to \mathrm{colim}_i W_i$$ are equal. Since  $\mathrm{colim}_i (W_i \times _{Y_i } W_i ) \to \mathrm{colim}_i W_i  \times _{\mathrm{colim}_i Y_i } \mathrm{colim}_i W_i $ is an epimorphism we conclude that the two legs of the pullback are equal.
  
Next we show stability of regular epis under pullbacks. Given a regular epi $f \colon X\to Z$ which occurs as the coequalizer of the pair of morphisms $h,\; k \colon W\to X$ and any morphsim $g \colon Y \to Z$ in $\mathcal{K}$ we are applying the lemma to the category $\mathcal{D}$ given as 
$$\xy \xymatrix{
& & y \ar[d]^{4} \\
w \ar@<1ex>[r]^{1} \ar@<-1ex>[r]_{2}&   x \ar[r]^{3}  & z}
\endxy$$
while $D\colon \mathcal{D} \to \mathcal{K}$ is defined by $Dw=W,$ $Dx=X,$ $Dz=Z,$ $D4=g,$ $D1=h,$ $D2 =k,$ $D3=f.$  It follows from the lemma that we can
write $X\cong \mathrm{colim}_i X_i ,$ $Y \cong \mathrm{colim}_i Y_i , $ $W \cong \mathrm{colim}_i W_i ,$ $h=\mathrm{colim}_i h_i ,$ $k= \mathrm{colim}_i k_i ,$ $f= \mathrm{colim}_i f_i$ and $g=\mathrm{colim}_i g_i$ over the same indexing category, with $X_i ,$ $Y_i ,$ $W_i$ in $\mathcal{C}.$ 
Then the coequalizers $q_ i \colon X_i \to Q_i$ of the pairs $(h_i , k_i )$ will have their codomains in $\mathcal{C}$
and they will factor as in the diagram

$$\xy \xymatrix{
W_i \ar@<1ex>[r]^{h_i} \ar@<-1ex>[r]_{k_i}&   X_i \ar[rr]^{f_i} \ar[dr]^{q_i} &  & Z_i \\
&  & Q_i \ar[ur]^{} &
} \endxy$$
Their pullbacks $X_i \times _{Z_i } Y_i  \twoheadrightarrow   Q_i \times _{Z_i } Y_i $ along the respective $Q_i \times _{Z_i} Y_i \to Q_i $ will be regular epimorphisms, hence the same will be 
$ \mathrm{colim}_i (X_i \times _{Z_i } Y_i ) \twoheadrightarrow   \mathrm{colim}_i (Q_i \times _{Z_i } Y_i ).$  
Using the commutation of colimits with coequalizers, the colimit $Q= \mathrm{colim}_i Q_i$ of these coequalizers will be 
$$\mathrm{colim}_i Q_i 
\cong \;   \mathrm{coeq} (h, k \colon W \to X )
\cong  Z$$
Hence we have a commutative diagram as below, with the lower rectangles being pullbacks. 

$$\xy \xymatrix{
 \mathrm{colim}_i (X_i \times _{Z_i } Y_i ) \ar@{->>}[d] \ar@{->>}[r]^{}  & \mathrm{colim}_i (Q_i \times _{Z_i } Y_i ) \ar@{->>}[d]^{}  & \\
X\times _Z Y  \ar[d]^{}  \ar[r]^{} &  Q\times _Z Y \ar[d]^{} \ar[r]^{\cong} &\mathrm{colim}_i Y_i  \ar[d]^{}\\
 \mathrm{colim}_i X_i  \ar[r]^{}  & \mathrm{colim}_i Q_i  \ar[r]^{\cong}&  \mathrm{colim}_i Z_i }
\endxy $$
The morphism  $X  \times _{Z }  Y  \to Q\times _Z Y $ is such that when composed with an epimorphism gives a regular epimorphism. Hence it is itself a regular epimorphism and so is its composition with the isomorphism $Q\times _Z Y \to Y.$ This proves the stability of regular epimorphisms under pullback.
\epf

Our intention is to apply the above to the category of compactly generated Hausdorff locales. We have seen that in that category, under the coreflectivity hypothesis, the colimit of the pullbacks of the compact sublocales of a compactly generated one along any morphism, maps epimorphically to the pullback of the locale along that morphism. We need epimorphicity of the comparison between the colimit of pullbacks of compact sublocales to the pullback of the colimits in order to use the above. To that end recall that when $\mathcal{I} \to \mathrm{Loc}$ is a directed system of inclusions of locales one has that the morphisms $X_i \to \mathrm{colim}_i X_i $ are also inclusions. More precisely, stated as a result about sup-lattices, following result appears in the proof of \cite{jt} Proposition I.2 and gives the corresponding result for frames.

\begin{lemma}
Let $(t_{ij} \colon A_i \to A_j )$ be an inverse directed diagram in the category of sup-lattices, such that all the transition maps $t_{ij}$ are surjective. Then the projections 
$$p_i \colon \mathrm{lim}_i A_i \to A_i $$ are also surjective. 
\end{lemma}

\begin{lemma}
Assume that a category $\mathcal{K}$ has the property that for a monomorphic directed diagram $X \colon \mathcal{I} \to \mathcal{K} $ over $Z$ and a morphism $Y \to Z,$ the canonical map 
$$ \mathrm{colim}_i ( X_i \times _Z Y )  \to (\mathrm{colim}_i X_i  ) \times _Z Y   $$ is a regular epimorphism. Then for a monomorphic directed system $X_i \rightarrow Z_i \leftarrow Y_i $ indexed by $\mathcal{I}$ the canonical comparison   
$$\mathrm{colim}_i  (X_i \times _{Z_i} Y_{i}  ) \rightarrow \mathrm{colim}_i X_i \times _{\mathrm{colim}_i Z_i} \mathrm{colim}_i Y_i $$
is a regular epimorphism. 
\end{lemma}

\pf 
Consider the pullback of the diagram
$$\xy \xymatrix{
 & \mathrm{colim}_i Y_i \ar[d]^{} \\
   \mathrm{colim}_i X_i  \ar[r]^{}  & \mathrm{colim}_i Z_i }
\endxy$$
Then 
$$\mathrm{colim}_i \mathrm{colim}_{i'} (X_i \times _{\mathrm{colim}_i Z_i} Y_{i'} ) 
 \; \cong \;  \mathrm{colim}_i  (X_i \times _{\mathrm{colim}_i Z_i} Y_{i}  ),$$
 by  directedness of $\mathcal{I}$ and we have regular epimorphisms
 $$\mathrm{colim}_i \mathrm{colim}_{i'} (X_i \times _{\mathrm{colim}_i Z_i} Y_{i'} ) \to \mathrm{colim}_i \; (X_i \times _{\mathrm{colim}_i Z_i} \mathrm{colim}_i Y_i ) \to \mathrm{colim}_i X_i \times _{\mathrm{colim}_i Z_i} \mathrm{colim}_i Y_i $$
by our assumption.  
Finally, since each $Z_i \to  \mathrm{colim}_i Z_i$ is monomorphism 
$$X_i \times _{\mathrm{colim}_i Z_i} Y_{i} \cong X_i \times _{Z_i } Y_{i} $$
as the following diagram of pullbacks indicates
$$\xy \xymatrix{
 X_i \times _{\mathrm{colim}_i Z_i} Y_{i} \ar[r]^{} \ar[d]^{} &  Y_i \ar[r]^{} \ar[d]^{} &  Y_i \ar[d]^{} \\
 X_i \ar[r]^{} \ar[d]^{} & Z_i \ar[r]^{id} \ar[d]^{id} & Z_i \ar[d]^{} \\
   X_i  \ar[r]^{} & Z_i \ar[r]^{} & \mathrm{colim}_i Z_i }
\endxy$$
\epf 

In view of the above Proposition and the previous lemma we get 

\begin{theorem}
Under the coreflectivity hypothesis, namely that it is coreflective in the category of Hausdorff locales, the category of compactly generated Hausdorff locales is regular. 
\end{theorem}

\pf
Apply Proposition 3.2 for $\mathcal{K}$ the category of compactly generated Hausdorff locales, $\mathcal{C}$ the category of compact Hausdorff locales, which is regular by \cite{t} 3.6.3. Recall that it is also closed under finite colimits inside the former. A finite coproduct of compact Hausdorff locales is obviously compact Hausdorff, while in the coequalizer  $q \colon X \to Q$ in $\mathrm{CGHLoc}$ of a pair of maps to a compact locale $X,$ $Q$ is the directed colimit of its compact sublocales. But $q[X]$ is one of them and it equals $Q.$
\epf 

\section{Effectivity of the category of compact Hausdorff locales}
Recall that a locale $X$ is regular if every element of its frame of opens is the supremum of all the elements of the frame that are well inside it. An element of a frame $U$ is well inside $V$, written $U \eqslantless V,$ if there exists a $W$ such that $U \wedge W= 0 $ and $W \vee V = X .$ Recall also that a locale is compact Hausdorff iff it is compact regular (a result due to \cite{v1}, see also \cite{t}, Theorem 3.4.2, for a different proof). A surjective map of locales is one where the inverse image of the corresponding map between the respective frames reflects order. 

\begin{proposition}
The image of a compact locale by a surjection is compact.
\end{proposition}

\pf 
Let $q \colon X \to Q$ be a surjection of locales, $q^* \colon OQ \to OX $ its inverse image and assume that $Q=\bigvee U_i ,$ where the union is directed. Then 
$$ X = q^* Q = q^* (\bigvee U_i ) = \bigvee q^* U_i $$
hence there is an $i$ such that $X=q^* U_i .$ It follows that 
$$Q = q_* q^* U_i = U_i ,$$ where the last equation follows by the fact that $q^*$ reflects  order.   
\epf

The following is Proposition 2 in \cite{hl}:

\begin{proposition}
The image of a regular locale by a proper surjection is regular.
\end{proposition}

\pf For a proper surjection $q \colon X \to Q$ with 
$X$ is regular we have that for every $V \in OQ$
$$q^* V = \bigvee \{U \in OX \; | \; U \eqslantless q^* V \}$$ from which we get 
$$V = q_* q^* V = q_* (\bigvee \{U \in OX \; | \; U \eqslantless q^* V \}) = \bigvee \{q_* U \in OX \; | \; U \eqslantless q^* V \}$$
since the involved supremum is directed hence preserved by $q_* $ (\cite{ss}, III 1.1).
Now $U \eqslantless q^* V$ implies $q_* U \eqslantless  V$ because if $W \in OX$ is a witness for the first relation, i.e we have 
$$ U \wedge W = 0 \qquad \mbox{and} \qquad q^* V \vee W = X$$ then 
$$ q_* U \wedge q_* W = q_* 0 =0 $$
(the latter because $Z\leq q_* 0 $ iff $q^* Z \leq 0 = q^* 0$ and $q^*$ reflects $\leq$)
and also 
$$Q = q_* X = q_* (q^* V \vee W) = V \vee q_* W$$
by properness of $q$. We conclude that $q_* W$ is a witness for $q_* U \eqslantless  V,$ hence 
$$V = \bigvee \{q_* U \in OQ \; | \; U \eqslantless q^* V \} \leq \bigvee \{q_* U \in OQ \; | \; q_* U \eqslantless  V \}.$$
\epf

\begin{theorem}
The category of compact Hausdorff locales is effective regular (= Barr-exact).
\end{theorem}

\pf First of all the category $\mathrm{CHausLoc}$ of compact Hausdorff locales is regular by \cite{t} 3.6.3. 
Equivalence relations in this category are proper and closed, as every map between compact Hausdorff locales is proper. 
We know from \cite{v} 5.17 that closed, proper equivalence relations on compact locales are effective, so they are the kernel pairs of their coequalizers in the category of locales. But the coequalizer of a proper equivalence relation is proper by \cite{v} 5.5. Hence the coequalizer in the category of locales of a (proper as it will be) equivalence relation between compact regular locales is compact regular, by the above two propositions. Since limits in the category in question are constructed as in the category of locales, we conclude that every equivalence relation in $\mathrm{CHausLoc},$ being proper, is the kernel pair of its coequalizer in the category of locales, which lives in $\mathrm{CHausLoc}.$
\epf

\begin{theorem}
The category of compact Hausdorff locales is a pretopos.
\end{theorem}

\pf We know from \cite{jt}, Proposition IV.4.1, that coproducts in the category of locales are universal. They are also disjoint since the pullback $X \times _{X + Y} Y$ of the injections of two locales into their coproduct is given by the frame $OX \bigotimes _{OX \times OY} OY$ that occurs as the tensor product in preframes of the corresponding frames over their product. Writing $p_1 ,$ $p_2$ for the projections of the product, an element $a \otimes b$ in the latter is then 
$$a\otimes b = (1\wedge p_1 (a,b)) \otimes b = 1 \otimes (p_2 (a,b) \wedge b ) = 1 \otimes b = 1 \otimes (p_2 (1,b) \wedge 1) = (1 \wedge p_1 (1,b)) \otimes 1 = 1 \otimes 1 $$
so the pullback is trivial.

Finite coproducts of compact Hausdorff locales are obviously compact Hausdorff themselves and, moreover, the category $\mathrm{CHausLoc}$ is closed in the category of locales under finite limits (in particular under pullbacks) \cite{t}, Lemma 3.6.3. Hence $\mathrm{CHausLoc}$ inherits from the category of locales the universality and disjointness of finite coproducts. \epf

\paragraph{Remark:} The category of compact Hausdorff spaces admits a characterization as the unique, up to equivalence, non-trivial, well-pointed, filtral pretopos with set-indexed copowers of its terminal object \cite{mr}. As the referee suggested, it would be interesting to know if the category of compact Hausdorff locales admits a similar ``pointless'' characterization.

\refs



\bibitem [Cagliari, Mantovani, Vitale 1995]{cmv} 
F.~Cagliari, S.~Mantovani, E.~M.~Vitale, Regularity of the category of Kelley spaces
{\em Applied Categorical Structures} 3 (1995) 357–361

\bibitem [Carboni, Pedicchio, Rosick\'{y} 2001]{cpr} A.~Carboni, M.~C.~Pedicchio, J.~Rosick\'{y}, Syntactic characterizations of various classes of locally presentable categories, 
{\em Journal of Pure and Applied Algebra} 161 (2001), 65-90

\bibitem [Chen 1991]{ch} X.~Chen, 
{\em Closed frame homomorphisms} PhD Thesis, McMaster University (1991)

\bibitem [Day and Street 1989]{ds} B.~Day and R.~Street, Localisation of locally presentable categories,
{\em Journal of Pure and Applied Algebra} 58 (1989), 227-233



\bibitem [Escardo 2006]{e} M.~Escardo, Compactly generated Hausdorff locales
{\em Annals of Pure and Applied Logic} 137 (2006), 147-163

\bibitem [He and Luo 2011]{hl} W.~He and M.~K.~Luo, A note on proper maps of locales
{\em Applied Categorical Structures} 19 (2011), 505-510

\bibitem [Johnstone 1982]{ss} P.~T.~Johnstone, {\em Stone Spaces}, Cambridge University Press (1982)

\bibitem [Johnstone 2002]{el} P.~T.~Johnstone, {\em Sketches of an Elephant: a topos theory  compendium:
vol.1 and vol.2}, Oxford Logic Guides 43 and 44. Clarendon Press,
Oxford (2002)

\bibitem [Joyal and Tierney 1984]{jt} A.~Joyal and M.~Tierney, An extension of the Galois theory of Grothendieck
{\em Memoirs of AMS} {\bf Vol. 51} (1984)



\bibitem [Marra and Reggio 2018]{mr} V.~Marra and L.~Reggio,
A characterisation of the category of compact Hausdorff spaces, arxiv.org 1808.09738

\bibitem [Townsend 1998]{t} C.~F.~Townsend, {\em Preframe techniques in constructive locale theory}, PhD Thesis, Imperial College (1998)

\bibitem [Vermeulen 1991]{v1} J.~J.~C.~Vermeulen, Some constructive results related to compactness and the (strong) Hausdorff property for locales, in
{\em Category Theory: Proceedings of the International Conference held in Como.} Lecture Notes in Mathematics 1488, Springer  (1991), 401-409

\bibitem [Vermeulen 1994]{v} J.~J.~C.~Vermeulen, Proper maps of locales
{\em Journal of Pure and Applied Algebra} 92 (1994), 79-107

\endrefs

\end{document}